\documentclass[fleqn,a4paper,11pt]{article}

\usepackage{amsmath,amssymb,amsthm}

\newtheorem{theorem}{Theorem}
\newtheorem{corollary}[theorem]{Corollary}
\newtheorem{proposition}[theorem]{Proposition}

\theoremstyle{remark}
\newtheorem{remark}{Remark}

\theoremstyle{definition}
\newtheorem{definition}{Definition}

\theoremstyle{definition}
\newtheorem{example}{Example}

\begin{document}

\title{A Remarkable Property of the Dynamic 
Optimization Extremals\footnote{Presented
at the contributed session \emph{Optimal Control and Calculus of Variations} 
of the 4th  International Optimization Conference in Portugal, \emph{Optimization 2001}, 
Aveiro, July 23--25, 2001. Accepted for publication in the journal 
\emph{Investiga\c{c}\~{a}o Operacional}, Vol. 22, Nr. 2, 2002, pp.~253--263.}}

\author{Delfim F. M. Torres\\
        \texttt{delfim@mat.ua.pt}}

\date{ R\&D Unit \emph{Mathematics and Applications}\\
       Department of Mathematics\\
       University of Aveiro\\
       3810-193 Aveiro, Portugal}

\maketitle


\begin{abstract}
We give conditions under which a function
$F\left(t,x,u,\psi_0,\psi\right)$ satisfies the relation
$\frac{\mathrm{d}F}{\mathrm{d}t} = \frac{\partial F}{\partial t}
+ \frac{\partial F}{\partial x} \cdot \frac{\partial H}{\partial \psi}
- \frac{\partial F}{\partial \psi} \cdot \frac{\partial H}{\partial x}$
along the Pontryagin extremals
$\left(x(\cdot),u(\cdot),\psi_0,\psi(\cdot)\right)$ of an optimal
control problem, where $H$ is the corresponding Hamiltonian. The
relation generalizes the well known fact that the equality
$\frac{\mathrm{d}H}{\mathrm{d}t} = \frac{\partial H}{\partial t}$
holds along the extremals of the problem, and that in the
autonomous case $H \equiv constant$. As applications of the new
relation, methods for obtaining conserved quantities along the 
Pontryagin extremals and for characterizing problems possessing
given constants of the motion are obtained.
\end{abstract}




\section{Introduction}

A dynamic optimization continuous problem poses the question of what
is the optimal magnitude of the choice variables, at each point of time,
in a given interval. To tackle such problems, three major approaches
are available: dynamic programming; the calculus of variations; and
the powerful and insightful optimal control. 
The calculus of variations is a classical subject, born in 1696 
with the brachistochrone problem, whose
field of applicability is broadened with 
optimal control theory. Dynamic programming is based
on the solution of a partial differential equation, 
known as the Hamilton-Jacobi-Bellman equation,
in order to compute a value function. Dynamic programming is
well designed to deal with optimization problems
in discrete time.
All these techniques are well known in the literature of
operations research (see e.g. \cite{bertsekas1,bertsekas2,valadares}), 
systems theory (see e.g. \cite{elgerd}),
economics (see e.g. \cite{chiang,leonard} and 
\cite[Cap\'{\i}tulo~14]{cesaltinaPires}) 
and management sciences (see e.g \cite{MalcolmGotterer})).
Here, we are concerned with the methods and procedures of optimal control.
This approach allows the effective study of
many optimization problems arising in such fields as
engineering, astronautics, mathematics, physics, economics,
business management and operations research,
due to its ability to deal with restrictions
on the variables and nonsmooth functions 
(see e.g. \cite{MalcolmGotterer,isaacs,flp,smith}).

At the core of optimal control
theory is the Pontryagin maximum principle -- the celebrated first
order necessary optimality condition -- whose solutions are called
(Pontryagin) extremals and which are obtained through a function $H$ called
Hamiltonian, akin to the Lagrangian function used in ordinary calculus
optimization problems (see e.g. \cite{pinch,smith})).
For autonomous problems of optimal control, i.e.
when the Hamiltonian $H$ does not depend explicitly
on time $t$, a basic property of the
Pontryagin extremals is the remarkable feature that the corresponding
Hamiltonian is constant along the extremals
(see e.g. \cite{pontryagin,gamkrelidze}). In classical
mechanics this property corresponds to energy conservation
(see e.g. \cite{lauwerier,rund}), while in the calculus of variations
it corresponds to the second Erdmann necessary optimality condition 
(see e.g. \cite{livroClarke}).
For problems of optimal control that depend upon time $t$ explicitly
(non-autonomous problems), the property amounts to the fact that
the total derivative with respect to time of the corresponding Hamiltonian
equals the partial derivative of the Hamiltonian with respect to
time:
\begin{eqnarray}
\label{eq:dHdt}
\frac{\mathrm{d}H}{\mathrm{d}t}\left(t,x(t),u(t),\psi_0,\psi(t)\right)
= \frac{\partial H}{\partial
t}\left(t,x(t),u(t),\psi_0,\psi(t)\right)
\end{eqnarray}
for almost all $t$ (see e.g. \cite{pontryagin,berkovitz,fattorini}).
This corresponds to the DuBois-Reymond necessary condition
of the calculus of variations (see e.g. \cite{cesari}).
Recent applications, in many different
contexts of the calculus of variations and optimal control, 
show the fundamental nature of the property \eqref{eq:dHdt}. 
It has been used in \cite{clarke85,ambrosio,sarychevTorres} to
establish Lipschitzian regularity of minimizers; in
\cite{clarke93} to establish some existence results; and
in \cite{3ncnw,ecc2001} to prove some generalizations of 
first Noether's theorem.
The techniques used in the proof of the relation are also
very useful, and have been applied in contexts
far away from dynamic optimization (see e.g. \cite{sarychev}).
In this note we give conditions under which a function
$F(t,x,u,\psi_0,\psi)$ satisfies the equality
\begin{eqnarray}
\label{eq:dFdt}
\frac{\mathrm{d}F}{\mathrm{d}t} = \frac{\partial F}{\partial t}
+ \frac{\partial F}{\partial x} \cdot \frac{\partial H}{\partial \psi}
- \frac{\partial F}{\partial \psi} \cdot \frac{\partial H}{\partial x} \, ,
\end{eqnarray}
almost everywhere, along the Pontryagin extremals. For $F = H$ equality (\ref{eq:dFdt})
reduces to (\ref{eq:dHdt}). 
As a corollary, we obtain a necessary and sufficient condition for
$F\left(t,x,u,\psi_0,\psi\right)$ to be a constant of the motion. 
From it, one is able to find constants of the motion that depend
on the control and that are not momentum maps, that is, one can
find preserved quantities $F\left(t,x(t),u(t),\psi_0,\psi(t)\right)$ 
along the Pontryagin extremals 
$\left(x(\cdot),u(\cdot),\psi_0,\psi(\cdot)\right)$ of the problem,
which are not of the form $\psi(t) \cdot C\left(x(t)\right)$.
This is in contrast with the results obtained in \cite{blankenstein},
where the conserved quantities are always
of the form $\psi(t) \cdot C\left(x(t)\right)$.
Our condition provides also a method for the characterization of optimal control 
problems with given constants of the motion. All these possibilities are
illustrated with examples.


\section{Preliminaries}

Without loss of generality (see e.g. \cite{berkovitz}), we will be
considering the optimal control problems in Lagrange form with fixed
initial time $a$ and fixed terminal time $b$ ($a < b$).


\subsection{Formulation of the Optimal Control Problem}

The problem consists of minimize a cost functional of the form
\begin{eqnarray}
\label{eq:J}
J\left[x(\cdot),u(\cdot)\right] = \int_a^b L\left(t,x(t),u(t)\right)
\mathrm{d}t \, ,
\end{eqnarray}
called the performance index,
among all the solutions of the vector differential equation
\begin{eqnarray}
\label{eq:dynamics}
\dot{x}(t) = \varphi\left(t,x(t),u(t)\right) \quad \text{ for a.a. }
t \in [a,b]\,  .
\end{eqnarray}
The \emph{state trajectory} $x(\cdot)$ is a $n$-vector absolutely continuous
function
\begin{eqnarray*}
x(\cdot) \in W_{1,1}\left([a,b];\mathbb{R}^n\right) \, ;
\end{eqnarray*}
and the \emph{control} $u(\cdot)$ is a $r$-vector measurable and bounded
function satisfying the control constraint $u(t) \in \Omega$,
\begin{eqnarray*}
u(\cdot) \in L_{\infty}\left([a,b];\Omega\right) \, .
\end{eqnarray*}
The set $\Omega \subseteq \mathbb{R}^r$ is called the \emph{control set}.
In general, the problem may include some boundary conditions and
state constrains, but they are not relevant for the present study:
the results obtained are independent of those restrictions. 
We assume the functions
$L : [a,b] \times \mathbb{R}^n \times \Omega \rightarrow \mathbb{R}$ and
$\varphi : [a,b] \times \mathbb{R}^n \times \Omega \rightarrow
\mathbb{R}^n$ to be continuous on $[a,b] \times \mathbb{R}^n \times \Omega$
and to have continuous derivatives with respect to $t$ and $x$.


\subsection{The Pontryagin Maximum Principle}

We shall now formulate the celebrated Pontryagin maximum principle
\cite{pontryagin}, which is a first-order necessary optimality
condition. The maximum principle provides a generalization
of the classical calculus of variations first-order necessary
optimality conditions and can treat problems in which upper and
lower bounds are imposed on the control variables -- a possibility
of considerable interest in operations research 
(see \cite{MalcolmGotterer}).

\begin{theorem}[Pontryagin maximum principle]
Let $\left(x(\cdot),u(\cdot)\right)$ be a minimizer of the optimal
control problem. Then, there exists a nonzero pair
$\left(\psi_0,\psi(\cdot)\right)$, where $\psi_0 \le 0$ is a constant
and $\psi(\cdot)$ a $n$-vector absolutely continuous function with
domain $[a,b]$, such that the following hold for almost all $t$ on
the interval $[a,b]$:
\begin{description}
\item[(i)] the Hamiltonian system
\begin{eqnarray*}
\left\{
\begin{array}{lcl}
\dot{x}(t) & = & \displaystyle
\frac{\partial H\left(t,x(t),u(t),\psi_0,\psi(t)\right)}{\partial
\psi} \, , \\[0.15in]
\dot{\psi}(t) & = & - \displaystyle
\frac{\partial H\left(t,x(t),u(t),\psi_0,\psi(t)\right)}{\partial x}
\, ;
\end{array}
\right.
\end{eqnarray*}
\item[(ii)] the maximality condition
\begin{eqnarray*}
H\left(t,x(t),u(t),\psi_0,\psi(t)\right) = \max_{v \in \Omega}
H\left(t,x(t),v,\psi_0,\psi(t)\right) \, ;
\end{eqnarray*}
\end{description}
with the Hamiltonian
$H(t,x,u,\psi_0,\psi) = \psi_0 L(t,x,u) + \psi \cdot \varphi(t,x,u)$.
\end{theorem}

\begin{definition}
A quadruple $\left(x(\cdot),\,u(\cdot),\,\psi_{0},\,\psi(\cdot)\right)$
satisfying the Hamiltonian system and
the maximality condition is called a (Pontryagin) extremal.
\end{definition}

\begin{remark}
Different terminology for the function $H$ can be found
in the literature. The Hamiltonian $H$ is sometimes called 
``unmaximized Hamiltonian'', ``pseudo-Hamiltonian'' or 
``Pontryagin function''.
\end{remark}

\begin{remark}
Transversality conditions may also appear in the Pontryagin maximum principle.
These conditions depend on the specific boundary conditions under consideration.
Our methods do not require the use of such transversality conditions and the 
results obtained are, as already mentioned, valid for arbitrary boundary conditions.
\end{remark}

\begin{remark}
The maximality condition is a static optimization problem.
The method of solving the optimal control problem 
\eqref{eq:J}--\eqref{eq:dynamics} via the maximum principle
consists of finding the solutions of the Hamiltonian system
by the elimination of the control with the aid of the maximality
condition. The required optimal solutions are found among these
extremals.
\end{remark}

The proof of the following theorem can be found,
for example, in \cite{pontryagin,berkovitz}.

\begin{theorem}
\label{r:teoremadHdt}
If $\left(x(\cdot),u(\cdot),\psi_0,\psi(\cdot)\right)$ is
a Pontryagin extremal, then the function 
$H\left(t,x(t),u(t),\psi_0,\psi(t)\right)$
is an absolutely continuous function of $t$ and satisfies
the equality (\ref{eq:dHdt}),
where on the left-hand side we have the total derivative with respect
to $t$, and on the right-hand side the partial derivative of the
Hamiltonian with respect to $t$.
\end{theorem}

As a particular case of Theorem~\ref{r:teoremadHdt}, 
when the Hamiltonian does
not depend explicitly on $t$, that is when the optimal control
problem is autonomous -- functions $L$ and $\varphi$ do not depend on
$t$ -- then the value of the Hamiltonian evaluated along an arbitrary
Pontryagin extremal $\left(x(\cdot),u(\cdot),\psi_0,\psi(\cdot)\right)$
of the problem turns out to be constant:
\begin{eqnarray*}
H(x(t),u(t),\psi_0,\psi(t)) \equiv const \, , \quad t \in [a,b]\, .
\end{eqnarray*}
We remark that Theorem~\ref{r:teoremadHdt} is a consequence of 
the Pontryagin maximum principle. We shall generalize Theorem~\ref{r:teoremadHdt}
in Section~\ref{s:mainresult}. Before, we review some facts from
functional analysis needed in the proof of our result.


\subsection{Facts from Functional Analysis}

First we introduce the concept of an absolutely continuous function
in $t$ uniformly with respect to $s$.

\begin{definition}
Let $\phi(s,t)$ be a real valued function defined on $[a,b] \times
[a,b]$. The function $\phi(s,t)$ is said to be an absolutely continuous
function in $t$ uniformly with respect to $s$ if, given $\varepsilon
> 0$, there exists $\delta > 0$, independent of $s$, such that for every
finite collection of disjoint intervals $(a_j,b_j) \subseteq [a,b]$
\begin{eqnarray*}
\sum_{j} \left(b_j - a_j\right) \le \delta \Rightarrow \sum_{j}
\left|\phi(s,b_j) - \phi(s,a_j)\right| \le \varepsilon
\quad \left(s \in [a,b]\right) \, .
\end{eqnarray*}
\end{definition}

The proof of the following two propositions can be found in
\cite[p. 74]{fattorini}.

\begin{proposition}
\label{proposition1}
Let $F\left(t,x,u,\psi_0,\psi\right)$,
$F : [a,b] \times \mathbb{R}^n \times \Omega \times \mathbb{R}_{0}^{-}
\times \mathbb{R}^n \rightarrow \mathbb{R}$, be continuously
differentiable with respect to $t$, $x$, $\psi$ for $u$ fixed, and
assume that there exists a function $G(\cdot) \in
L_1\left([a,b] ; \mathbb{R} \right)$ such that
\begin{eqnarray*}
\left\|\nabla_{(t,x,\psi)}
F\left(t,x(t),u(s),\psi_0,\psi(t)\right)\right\| \le G(t) \, 
\quad \left(s, t \in [a,b]\right) \, .
\end{eqnarray*}
Then $\phi(s,t) = F\left(t,x(t),u(s),\psi_0,\psi(t)\right)$ is
absolutely continuous in $t$ uniformly with respect to $s$ on
$[a,b]$.
\end{proposition}

\begin{proposition}
\label{proposition2}
Let $\phi(s,t)$, $\phi : [a,b] \times [a,b] \rightarrow \mathbb{R}$,
be an absolutely continuous function in $t$ uniformly with respect to
$s$ satisfying
\begin{eqnarray*}
\phi(t,t) = \max_{s \in [a,b]} \phi(s,t)
\end{eqnarray*}
in a set dense in $[a,b]$. Then the function $\phi(t,t)$ can be
uniquely extended to a function $m(t)$ absolutely continuous on
$[a,b]$.
\end{proposition}


\section{Main Result}
\label{s:mainresult}

Our result is a generalization of the 
Theorem~\ref{r:teoremadHdt}.

\begin{theorem}
\label{r:mainresult}
If $F\left(t,x,u,\psi_0,\psi\right)$ is a real valued function as in
Proposition~\ref{proposition1} and besides satisfies
\begin{eqnarray}
\label{eq:MAX}
F\left(t,x(t),u(t),\psi_0,\psi(t)\right)
= \max_{v \in \Omega} F\left(t,x(t),v,\psi_0,\psi(t)\right)
\end{eqnarray}
a.e. in $t \in [a,b]$ along any Pontryagin extremal
$\left(x(\cdot),u(\cdot),\psi_0,\psi(\cdot)\right)$ of the optimal
control problem, then
$t \rightarrow F\left(t,x(t),u(t),\psi_0,\psi(t)\right)$ is
absolutely continuous and the equality
\begin{eqnarray}
\label{eq:MYEQ}
\frac{\mathrm{d}F}{\mathrm{d}t} = \frac{\partial F}{\partial t}
+ \frac{\partial F}{\partial x} \cdot \frac{\partial H}{\partial \psi}
- \frac{\partial F}{\partial \psi} \cdot \frac{\partial H}{\partial x}
\end{eqnarray}
holds along the extremals.
\end{theorem}

\begin{proof}
Our proof is an extension of the standard proof of 
Theorem~\ref{r:teoremadHdt}.
Let $\left(x(\cdot),u(\cdot),\psi_0,\psi(\cdot)\right)$ be a
Pontryagin extremal of the problem. Setting $v = u(s)$ in
(\ref{eq:MAX}) we obtain that
$\phi(s,t) = F\left(t,x(t),u(s),\psi_0,\psi(t)\right)$ satisfies
\begin{eqnarray}
\label{eq:phittgephist}
\phi(t,t) \ge \phi(s,t) \, , \quad s \in [a,b] \, ,
\end{eqnarray}
for $t$ in a set of full measure on $[a,b]$. Proposition~\ref{proposition2}
then implies that $m(t) = \phi(t,t)
= F\left(t,x(t),u(t),\psi_0,\psi(t)\right)$ is an absolutely
continuous function on $[a,b]$. It remains to prove that
\begin{eqnarray*}
\dot{m}(t) = \frac{\partial F}{\partial t}\left(\pi(t)\right)
+ \frac{\partial F}{\partial x}\left(\pi(t)\right) \cdot
\frac{\partial H}{\partial \psi}\left(\pi(t)\right)
- \frac{\partial F}{\partial \psi}\left(\pi(t)\right) \cdot
\frac{\partial H}{\partial x}\left(\pi(t)\right) \, ,
\end{eqnarray*}
where $\pi(t) = \left(t,x(t),u(t),\psi_0,\psi(t)\right)$. Since
\begin{equation*}
\frac{m(t+h)-m(t)}{h} = \frac{\phi(t+h,t+h) - \phi(t,t+h)}{h} +
                          \frac{\phi(t,t+h) - \phi(t,t)}{h}
\end{equation*}
and by the hypotheses the left-hand side and the second term on the
right-hand side have a limit as $h \rightarrow 0$, one concludes that
the first term on the right must have a limit as well. From
(\ref{eq:phittgephist}) $\phi(t+h,t+h) \ge \phi(t,t+h)$ and it
follows that $\frac{\phi(t+h,t+h) - \phi(t,t+h)}{h}$ is nonnegative
when $h > 0$ and nonpositive when $h < 0$; thus, its limit must be
zero when $h \rightarrow 0$. In this way we obtain that
\begin{eqnarray*}
\dot{m}(t) &=& \lim_{h \rightarrow 0}
\frac{F\left(t+h,x(t+h),u(t),\psi_0,\psi(t+h)\right) -
F\left(t,x(t),u(t),\psi_0,\psi(t)\right)}{h} \\
&=& \frac{\partial F}{\partial t}\left(\pi(t)\right)
+ \frac{\partial F}{\partial x}\left(\pi(t)\right) \cdot \dot{x}(t)
+ \frac{\partial F}{\partial \psi}\left(\pi(t)\right) \cdot
\dot{\psi}(t) \, ,
\end{eqnarray*}
and the conclusion follows from the Hamiltonian system.
\end{proof}

\begin{corollary}
Let $F\left(t,x,u,\psi_0,\psi\right)$,
$F : [a,b] \times \mathbb{R}^n \times \Omega \times \mathbb{R}_{0}^{-}
\times \mathbb{R}^n \rightarrow \mathbb{R}$, be continuously
differentiable with respect to $t$, $x$, $\psi$ for $u$ fixed;
and $\left(x(\cdot),u(\cdot),\psi_0,\psi(\cdot)\right)$ be an
extremal. If
\begin{description}
\item[(i)] $F\left(t,x(t),u(t),\psi_0,\psi(t)\right)$ is
absolutely continuous in $t$;
\item[(ii)] $F\left(t,x(t),u(t),\psi_0,\psi(t)\right)
= \displaystyle \max_{v \in \Omega} F\left(t,x(t),v,\psi_0,\psi(t)\right)$
a.e. in $a \le t \le b$; 
\end{description}
then the equality \eqref{eq:MYEQ} holds along the extremal.
\end{corollary}

Possible applications of Theorem~\ref{r:mainresult} 
follow in the next section.


\section{Applications of the Main Result}

Solving the Hamiltonian system by 
the elimination of the control with the aid of the maximality condition
is typically a difficult task. Therefore, it is worthwhile to look
for circumstances which make the solution easier. This is the case
when the extremals don't change the value of a given function. 
Indeed, the existence of such a function, called constant of the motion, 
may be used for reducing the dimension of the Hamiltonian system
(see e.g. \cite[M\'{o}dulo 5]{vasile}).
In extreme cases, with a sufficiently large number of (independent) 
constants of the motion, one can solve the problem completely.


\subsection{Constants of the Motion}

From Theorem~\ref{r:mainresult}, one immediately obtains a necessary
and sufficient condition for a function to be a constant of the motion.

\begin{definition}
A quantity $F(t,x,u,\psi_0,\psi)$ which is constant along every
Pontryagin extremal $\left(x(\cdot),u(\cdot),\psi_0,\psi(\cdot)\right)$ of
the problem, is called a \emph{constant of the motion}.
\end{definition}

\begin{corollary}
\label{r:corollary:CM}
Under the conditions of Theorem~\ref{r:mainresult},
$F(t,x,u,\psi_0,\psi)$ is a constant of the motion if and only if
\begin{equation}
\label{e:necsufCM}
\frac{\partial F}{\partial t}
+ \frac{\partial F}{\partial x} \cdot \frac{\partial H}{\partial \psi}
- \frac{\partial F}{\partial \psi} \cdot \frac{\partial H}{\partial x}
= 0
\end{equation}
holds, almost everywhere, along the Pontryagin extremals of the optimal control problem.
\end{corollary}

\begin{example} ($n = 4$, $r = 2$, $\Omega = \mathbb{R}^2$)
Let us consider the problem
\begin{gather*}
\int_a^b \left(\left(u_1(t)\right)^2 + \left(u_2(t)\right)^2 \right) \, dt 
\longrightarrow \min \, , \\
\begin{cases}
\dot{x_1}(t) = x_3(t) \\
\dot{x_2}(t) = x_4(t) \\
\dot{x_3}(t) = - x_1(t) \left(\left(x_1(t)\right)^2 + 
\left(x_2(t)\right)^2\right) + u_1(t) \\
\dot{x_4}(t) = - x_2(t) \left(\left(x_1(t)\right)^2 + 
\left(x_2(t)\right)^2\right) + u_2(t) \, .
\end{cases}
\end{gather*}
The corresponding Hamiltonian function is
\begin{multline*}
H\left(x_1,x_2,x_3,x_4,u_1,u_2,\psi_0,\psi_1,\psi_2,\psi_3,\psi_4\right) 
= \psi_0 \left(u_{1}^2+u_{2}^2\right) + \psi_1 x_3 \\
+ \psi_2 x_4 
- \psi_3 x_1 \left(x_{1}^2+x_{2}^2\right) + \psi_3 u_1
- \psi_4 x_2 \left(x_{1}^2+x_{2}^2\right) + \psi_4 u_2 \, .
\end{multline*}
We claim that
\begin{equation}
\label{e:CMex1}
F = - \psi_1 x_2 + \psi_2 x_1 - \psi_3 x_4 + \psi_4 x_3
\end{equation}
is a constant of the motion for the problem. Direct
calculations show that
\begin{equation}
\label{e:calculosEX1}
\frac{\partial F}{\partial t} + 
\sum_{i=1}^{4} \frac{\partial F}{\partial x_i} \frac{\partial H}{\partial \psi_i}
- \sum_{i=1}^{4} \frac{\partial F}{\partial \psi_i} \frac{\partial H}{\partial x_i} =
\psi_4 u_1 - \psi_3 u_2 \, .
\end{equation}
From the maximality condition it follows that $\frac{\partial H}{\partial u_1} = 0$
and $\frac{\partial H}{\partial u_2} = 0$, that is, $2 \psi_0 u_1 + \psi_3 = 0$
and $2 \psi_0 u_2 + \psi_4 = 0$. Using these last two identities in \eqref{e:calculosEX1}
one concludes from Corollary~\ref{r:corollary:CM} that \eqref{e:CMex1} is a
constant of the motion.
\end{example}


\subsection{Characterization of Optimal Control Problems }

We shall endeavor here to find a method to synthesize optimal control 
problems with given constants of the motion.
If a function $F$ is fixed \emph{a priori}, we can regard equality
\eqref{e:necsufCM} as a partial differential equation in the unknown
Hamiltonian $H$. Obviously, if this differential equation admits a
solution, then an optimal control problem can be constructed with the
constant of the motion $F$. We shall illustrate the general idea in
special situations.

\begin{example}
The Hamiltonian $H$ is a constant of the motion if and only if
$\frac{\partial H}{\partial t} = 0$. Condition is trivially 
satisfied for autonomous problems.
\end{example}

\begin{example}
Function $\psi x + H t$ is a constant of the motion if and only if
$H = \frac{\partial H}{\partial x} x 
- \frac{\partial H}{\partial \psi} \psi
- \frac{\partial H}{\partial t} t$. Condition is satisfied,
for example, for problems of the form ($0 < a < b$)
\begin{gather*}
\int_a^b \frac{L\left(t x(t),u(t)\right)}{t} dt \longrightarrow \min \, , \\
\dot{x}(t) = \frac{\varphi\left(t x(t),u(t)\right)}{t^2} \, .
\end{gather*}
\end{example}

\begin{example}
We conclude from Corollary~\ref{r:corollary:CM} that
a necessary and sufficient condition for $H \psi x$ to be a 
constant of the motion is
\begin{equation*}
\psi x \frac{\partial H}{\partial t}
+ \psi H \frac{\partial H}{\partial \psi}
- H x \frac{\partial H}{\partial x} = 0 \, .
\end{equation*}
A simple problem with constant of the motion $H \psi x$ is therefore
\begin{gather*}
\int_a^b L\left(u(t)\right) dt \longrightarrow \min \, , \\
\dot{x}(t) = \varphi\left(u(t)\right) x(t) \, .
\end{gather*}
\end{example}

\begin{example}
The following optimization problem is important in the
study of cubic polynomials on Riemannian manifolds
(see \cite[p. 39]{camarinha} and \cite{LeiteCamarinhaCrouch}).
Here we consider the particular case when one has
2-dimensional state and $n$ controls:
\begin{gather}
\int_0^T \left(\left(u_1(t)\right)^2 + \cdots + \left(u_n(t)\right)^2 \right) dt 
\longrightarrow \min \, , \label{e:RM-SC} \\
\begin{cases}
\dot{x_1}(t) = x_2(t) \, , \\
\dot{x_2}(t) = X_1\left(x_1(t)\right) u_1(t) + \cdots + X_n\left(x_1(t)\right) u_n(t) \, .
\end{cases} \notag
\end{gather}
Functions $X_i(\cdot)$, $i=1,\,\ldots,\,n$, are assumed smooth. 
The Hamiltonian for the problem is
\begin{equation*}
H = \psi_0 \left(u_1^2 + \cdots + u_n^2\right) + \psi_1 x_2 +
\psi_2 \left(X_1(x_1) u_1 + \cdots + X_n(x_1) u_n\right) \, .
\end{equation*}
As far as the problem is autonomous, the Hamiltonian is a constant of the
motion. We are interested in finding a new constant of the motion
for the problem. We will look for one of the form
\begin{equation*}
F = k_1 \psi_1 x_1 + k_2 \psi_2 x_2 \, ,
\end{equation*}
where $k_1$ and $k_2$ are constants. This is a typical constant of
the motion, known in the literature by \emph{momentum map} 
(see \cite{blankenstein}). First we note that
\begin{equation*}
\frac{\partial F}{\partial t} = 0 \, , 
\frac{\partial F}{\partial x_1} = k_1 \psi_1 \, , 
\frac{\partial F}{\partial x_2} = k_2 \psi_2 \, , 
\frac{\partial F}{\partial \psi_1} = k_1 x_1 \, , 
\frac{\partial F}{\partial \psi_2} = k_2 x_2 \, ,
\end{equation*}
and
\begin{gather*}
\frac{\partial H}{\partial x_1} = \psi_2 
\left(X'_1(x_1) u_1 + \cdots + X'_n(x_1) u_n\right) \, , 
\frac{\partial H}{\partial x_2} = \psi_1 \, , \\
\frac{\partial H}{\partial \psi_1} = x_2  \, ,
\frac{\partial H}{\partial \psi_2} = X_1(x_1) u_1 + \cdots + X_n(x_1) u_n \, .
\end{gather*}
Substituting these quantities into \eqref{e:necsufCM}
we obtain that
\begin{multline*}
k_1 \psi_1 x_2 + k_2 \psi_2 \left(X_1(x_1) u_1 + \cdots + X_n(x_1) u_n\right) \\
- k_1 x_1 \psi_2 \left(X'_1(x_1) u_1 + \cdots + X'_n(x_1) u_n\right)
- k_2 x_2 \psi_1 = 0 \, .
\end{multline*}
The equality is trivially satisfied if $k_1 = k_2$ and
$X'_i(x_1) x_1 = X_i(x_1)$, $i = 1,\,\ldots,\,n$.
We have just proved the following proposition.

\begin{proposition}
If the homogeneity condition $X_i\left(\lambda x_1\right) = \lambda X_i(x_1)$
($i = 1,\,\ldots,\,n$), $\forall \, \lambda > 0$, holds, then
$\psi_1(t) x_1(t) + \psi_2(t) x_2(t)$
is constant in $t \in [0,T]$ along the extremals of the
problem \eqref{e:RM-SC}.
\end{proposition}
\end{example}


\section*{Acknowledgments}

The author is in debt to A.~V.~Sarychev
for the many useful advises, comments and suggestions.
The research was supported by
the program PRODEP III 5.3/C/200.009/2000.




\begin{thebibliography}{99}

\bibitem{ambrosio} Ambrosio L., Ascenzi O., Buttazzo G.
Lipschitz Regularity for Minimizers of Integral Functionals
with Highly Discontinuous Integrands.
J. Math. Anal. Appl. 142, 1989, pp.~301--316.

\bibitem{berkovitz} Berkovitz L.~D. 
Optimal Control Theory.
Applied Mathematical Sciences 12,
Springer-Verlag, New York, 1974.

\bibitem{bertsekas1} Bertsekas D.~P. 
Dynamic Programming and Optimal Control, Vol.~I
(2nd ed.). Athena Scientific, Belmont, Massachusetts, 2000.

\bibitem{bertsekas2} Bertsekas D.~P. 
Dynamic Programming and Optimal Control, Vol.~II.
Athena Scientific, Belmont, Massachusetts, 1995.

\bibitem{blankenstein} Blankenstein~G., van der Schaft~A.
Optimal control and implicit Hamiltonian systems. In:
Isidori~A., Lamnabhi-Lagarrigue~F., Respondek~W. (eds).
Nonlinear control in the year 2000, vol.~1 (Paris). 
Springer, London. 2001, pp.~185--205.

\bibitem{camarinha} Camarinha M. 
A Geometria dos Polin\'{o}mios C\'{u}bicos
em Variedades Riemannianas. Ph.D. thesis,
Departamento de Matem\'{a}tica, Universidade de Coimbra, 
Coimbra, 1996.

\bibitem{cesari} Cesari L.
Optimization---Theory and Applications.
Springer-Verlag, New York, 1983.

\bibitem{chiang} Chiang A.~C. 
Elements of Dynamic Optimization.
McGraw-Hill Inc, 1992.

\bibitem{livroClarke} Clarke F.~H. 
Optimization and Nonsmooth Analysis.
John Wiley \& Sons Inc., New York, 1983.

\bibitem{clarke93} Clarke F.~H. 
An Indirect Method in the Calculus of Variations.
Trans. Amer. Math. Soc. 336, 1993, pp.~655--673.

\bibitem{clarke85} Clarke F.~H., Vinter R.~B. 
Regularity Properties of Solutions to the Basic Problem in the
Calculus of Variations. Trans. Amer. Math. Soc. 
289, 1985, pp.~73--98.

\bibitem{MalcolmGotterer} Connors M.~M., Teichroew D. 
Optimal Control of Dynamic Operations Research Models.
International Textbook Company, Scranton, Pennsylvania, 1967.

\bibitem{elgerd} Elgerd O.~I. 
Control Systems Theory. McGraw-Hill Inc, 1967.

\bibitem{fattorini} Fattorini H.~O. 
Infinite Dimensional Optimization and Control Theory.
Encyclopedia of Mathematics and Its Applications 62,
Cambridge University Press, Cambridge, 1999.

\bibitem{sarychev} Freiling G., Jank G., Sarychev A. 
Non-blow-up Conditions for Riccati-type Matrix Differential
and Difference Equations. Results Math. 37, 2000, pp.~84--103.

\bibitem{gamkrelidze} Gamkrelidze R.~V. 
Principles of Optimal Control Theory.
Mathematical Concepts and Methods in Science and Engineering
7, Plenum Press, New York, 1978.

\bibitem{isaacs} Isaacs R. 
Differential Games -- A Mathematical Theory with Applications
to Warfare and Pursuit, Control and Optimization.
Dover Publications Inc., Mineola, New York, 1999.

\bibitem{lauwerier} Lauwerier H.~A. 
Calculus of Variations in Mathematical Physics.
Mathematical Centre Tracts 14, Mathematisch Centrum, 
Amsterdam, 1966.

\bibitem{leonard} L\'{e}onard D., Van Long N. 
Optimal Control Theory and Static Optimization in Economics.
Cambridge University Press, Cambridge, 1992.

\bibitem{flp} Pereira F.~L.
Control Design for Autonomous Vehicles:
A Dynamic Optimization Perspective.
European Journal of Control 7, 2001, pp.~178--202.

\bibitem{pinch} Pinch E.~R. 
Optimal Control and the Calculus of Variations.
Oxford University Press, Oxford, 1995.

\bibitem{cesaltinaPires} Pires C.
C\'{a}lculo para Economistas.
McGraw-Hill de Portugal Lda., 2001.

\bibitem{pontryagin} Pontryagin L.~S., Boltyanskii V.~G.,
Gamkrelidze R.~V., Mischenko E.~F. 
The Mathematical Theory of Optimal Processes. 
John Wiley, New York, 1962.

\bibitem{rund} Rund H. 
The Hamilton--Jacobi Theory in the Calculus of Variations,
Its Role in Mathematics and Physics.
D.~Van Nostrand Co., Ltd., London--Toronto, 
Ont.--New York, 1966.

\bibitem{sarychevTorres} Sarychev A.~V., Torres D.~F.~M.
Lipschitzian Regularity of Minimizers for Optimal Control
Problems with Control-Affine Dynamics.
Applied Mathematics and Optimization, 41, 2000,
pp.~237--254.

\bibitem{LeiteCamarinhaCrouch} Silva Leite F., Camarinha M., Crouch P.
Elastic Curves as Solutions of Riemannian and Sub-Riemannian Control Problems.
Math. Control Signals Systems 13, 2000, pp.~140--155.

\bibitem{smith} Smith D.~R. 
Variational Methods in Optimization.
Dover Publications Inc., Mineola, New York, 1998.

\bibitem{vasile} Staicu V.
Equa\c{c}\~{o}es Diferenciais. Relat\'{o}rio da disciplina
de Equa\c{c}\~{o}es Diferenciais, Provas de Agrega\c{c}\~{a}o em
Matem\'{a}tica, Universidade de Aveiro, 2000.

\bibitem{3ncnw} Torres D.~F.~M. 
Conservation Laws in Optimal Control. 
Dynamics, Bifurcations and Control, 
Lecture Notes in Control and Information Sciences 273, 
Springer-Verlag, Berlin, Heidelberg, 2002, pp.~287--296.

\bibitem{ecc2001} Torres D.~F.~M. 
On the Noether Theorem for Optimal Control.
European Journal of Control, 8(1) 2002, pp.~56--63. 

\bibitem{valadares} Valadares Tavares L., Nunes Correia F. 
Optimiza\c{c}\~{a}o Linear e N\~{a}o Linear -- Conceitos, M\'{e}todos
e Algoritmos. Funda\c{c}\~{a}o Calouste Gulbenkian, Lisboa, 1986.

\end{thebibliography}
\end{document}